# The maximum variance of a finite dataset, given its mean, minimum, and maximum

Jules L. Ellis[a]

[a]*Department of Theory, Methods and Statistics, Open Universiteit, Heerlen, The Netherlands*

**Author note**

Jules L. Ellis https://orcid.org/0000-0002-4429-5475

E-mail: jules.ellis@ou.nl

**Funding statement.** This research received no specific grant funding form any funding agency, commercial or not-for-profit sectors.

**Competing interests.** The author declares none.

**MSC 2020-codes.** 60E15 — Inequalities; stochastic orderings; 26D15 — Inequalities for sums, series and integrals; 26D20 — Other analytical inequalities

**The maximum variance of a finite dataset, given its mean, minimum, and maximum**

This paper derives the maximum variance of a finite dataset of real numbers, given their mean, minimum and maximum. An example is provided in which the maximum variance is less than half of the Bhatia-Davis upper bound, (maximum - mean)(mean - minimum). As the dataset length increases, the maximum variance under these constraints approaches this bound from below.

**Introduction**

Consider a finite dataset (or sequence) of real numbers with known length, lower and upper bound, and (arithmetic) mean. What is the maximum value of the variance of the dataset under these constraints? Bhatia and Davis [1] showed that the variance is at most (maximum – mean)(mean – minimum). Other authors [2–5] have provided related inequalities, but, to my knowledge, none are claimed to be sharp upper bounds on the variance for the problem under discussion. I will provide a sharp upper bound. For example, suppose we have five numbers between 0 and 1 with mean 0.1. What is the maximum possible variance of these five numbers? The Bhatia-Davis upper bound is $0.1 * (1 - 0.1) = 0.09$, but the actual maximum is 0.04, as will be discussed below.

For any dataset $x = (x_1, \ldots, x_n) \in \mathbb{R}^n$, denote the mean as $\bar{x} = \frac{\sum_{i=1}^n x_i}{n}$ and let the variance be defined as

$$\mathrm{var}(x) = \frac{\sum_{i=1}^n (x_i - \bar{x})^2}{n} = \frac{\sum_{i=1}^n x_i^2}{n} - \bar{x}^2.$$

Suppose that it is known that the $x_i$ are bounded by $m \leq x_i \leq M$, and that their mean is known to be equal to $c$. The question is: what is the maximum of $\mathrm{var}(x)$ under the given constraints? For ease of expression, we will study the problem first with $m = 0$ and $M = 1$.

## Which datasets maximize the variance under the given constraints?

Without the constraint on the mean, the variance is maximal only if all $x_i$ are equal to one of the boundary values. For example, if $0 \leq x_i \leq 1$ for $i = 1, \ldots, n$, the variance is maximal only if all $x_i$ are either 0 or 1. This condition, however, implies that $\bar{x}$ is a multiple of $1/n$, and therefore it cannot be a solution that maximizes the variance under the constraint that, say $\bar{x} = 0.5$ when $n = 5$.

**Lemma 1.** *Let $\boldsymbol{x}$ be the dataset $\boldsymbol{y}$ that maximizes $\mathrm{var}(\boldsymbol{y})$ under the constraints $\boldsymbol{y} \in [0, 1]^n$ and $\bar{y} = c$. Then there is at most one $x_i$ with $0 < x_i < 1$.*

*Proof.* By contradiction. Suppose there is more than one $x_i$ with $0 < x_i < 1$. Without loss of generality we may assume that these are $x_1$ and $x_2$, with $x_1 \leq x_2$. Let $\varepsilon = \min(x_1, 1 - x_2)$, and define a vector $\boldsymbol{y} = (y_1, \ldots, y_n)$ as follows:

$$y_1 = x_1 - \varepsilon$$

$$y_2 = x_2 + \varepsilon$$

$$y_i = x_i \text{ for } i = 3, \ldots, n.$$

Then $\boldsymbol{y} \in [0, 1]^n$ and $\frac{\sum_{i=1}^{n} y_i}{n} = c$, and

$$\mathrm{var}(\boldsymbol{y}) = \frac{\sum_{i=1}^{n} y_i^2}{n} - c^2 =$$

$$\frac{y_1^2 + y_2^2 + \sum_{i=3}^{n} y_i^2}{n} - c^2 =$$

$$\frac{(x_1 - \varepsilon)^2 + (x_2 + \varepsilon)^2 + \sum_{i=3}^{n} x_i^2}{n} - c^2 =$$

$$\frac{2\varepsilon^2 + 2(x_2 - x_1)\varepsilon}{n} + \frac{\sum_{i=1}^{n} x_i^2}{n} - c^2.$$

But $\varepsilon > 0$ and $x_2 - x_1 \geq 0$; therefore $\text{var}(y) > \text{var}(x)$, contradicting the premise that $\text{var}(x)$ is maximal. ∎

**Lemma 2.** *Suppose $x$ is the dataset $y$ that maximizes $\text{var}(y)$ under the constraints $y \in [0,1]^n$ and $\bar{y} = c$. Let $a$ be the fractional part of $nc$. If there is an $x_i$ with $0 < x_i < 1$, then $x_i = a$.*

*Proof.* According to Lemma 1, there can be at most one $x_i$ that is not equal to 0 or 1. Denote this $x_i$ as $x_s$. Let the number of $x_i$ with $x_i = 1$ be denoted as $k$. Now $k + x_s = nc$. Therefore, $x_s$ must be equal to the fractional part of $nc$, i.e., $x_s = a$. ∎

**Example 1.** Consider the maximum variance of $x_1, \ldots, x_5$ under the constraints $0 \leq x_1, \ldots, x_5 \leq 1$ and $\frac{\sum_{i=1}^{5} x_i}{5} = 0.1$. Using Lemma 1, we know that four out of five $x_i$ are equal to 0 or 1. This must be four 0s, otherwise the mean cannot be less than 0.2. The only non-zero element must then be 0.5. Therefore, we find the maximum variance under the given constraints as the variance of the dataset $(0, 0, 0, 0, 0.5)$, which is 0.04. This is less than the bound of Bhatia and Davis [1], which is $0.1 * (1 - 0.1) = 0.09$. This idea to find the maximum is generalized in the next section.

**The maximum of the variance under the given constraints**

**Theorem 1.** *Let $x$ be the dataset $y$ that maximizes $\text{var}(y)$ under the constraints $y \in [0,1]^n$ and $\bar{y} = c$, and let $a$ be the fractional part of $nc$. Then*

$$\text{var}(x) = c(1-c) - \frac{1}{n} a(1-a).$$

*Proof.* According to Lemma 1, there can be at most one $x_i$ that is not equal to 0 or 1. If there are no such $x_i$, then all $x_i$ are equal to either 0 or 1, which implies that $nc$ is integer and $a = 0$, and that $\text{var}(x) = c(1-c)$; in agreement with the conclusion of the theorem. In the rest of the proof, we consider the case where there is one $x_i$ with $0 < x_i < 1$. Denote this $x_i$ as $x_s$. Let the number of $x_i$ with $x_i = 1$ be denoted as $k$. Now $k + x_s = nc$, and by Lemma 2, $x_s = a$, so

$$\frac{k}{n} = c - \frac{a}{n}.$$

The variance can thus be obtained as

$$\text{var}(x) = \frac{\sum_{i=1}^{n}(x_i - c)^2}{n}$$

(substitute $k$ values 1 and $n - k - 1$ values 0 and one value $a$ for for the $x_i$:)

$$= \frac{k}{n}(1 - c)^2 + \frac{n - k - 1}{n}c^2 + \frac{1}{n}(a - c)^2$$

$$= \frac{k}{n} - \frac{k}{n}2c + c^2 + \frac{1}{n}a^2 - \frac{1}{n}2ca$$

(substitute $\frac{k}{n} = c - \frac{a}{n}$:)

$$= c - \frac{a}{n} - (c - \frac{a}{n})2c + c^2 + \frac{1}{n}a^2 - \frac{1}{n}2ca$$

$$= c(1 - c) - \frac{1}{n}a(1 - a). \blacksquare$$

Since $0 \leq a < 1$, we have $0 \leq a(1 - a) \leq \frac{1}{4}$. Therefore, we have the following Corollary.

**Corollary 1.** *Let $x$ be the dataset $y$ that maximizes $\text{var}(y)$ under the constraints $y \in [0, 1]^n$ and $\bar{y} = c$. Then*

$$c(1 - c) - \frac{1}{4n} \leq \text{var}(x) \leq c(1 - c).$$

**Example 2**. If $x \in [0,1]^n$ then Theorem 1 implies that $\sum_{i=1}^{n} x_i^2 \leq \sum_{i=1}^{n} x_i - a(1-a)$, where $a$ is the fractional part of $\sum_{i=1}^{n} x_i$.

Finally, consider the case where the dataset length $n$ is unknown and may be arbitrarily large, and where the bounds are $m$ and $M$ rather than 0 and 1.

**Corollary 2**. *As* $n \to \infty$,

$$\max_{x \in [m,M]^n}\{\text{var}(x) : \bar{x} = c\} \uparrow (M-c)(c-m).$$

This follows directly from Theorem 1 with the transformation $x \mapsto (M-m)x + m$. Thus, given a finite dataset of real numbers with known bounds and average, its variance can be arbitrarily close to the upper bound (maximum – mean)(mean – minimum) if the dataset is long enough.

**Example 3**. The coefficient of variation (CV) or relative standard deviation is defined as $\sqrt{\text{var}(x)}/\bar{x}$. For a finite dataset of numbers between 0 and 1, Theorem 1 implies that the CV is at most $\sqrt{\bar{x}^{-1} - 1}$, and this is the tightest bound that is possible without knowing the length of the dataset.